
\input amstex
\documentstyle{amsppt}
\def\g{\gamma}
\def\k{\kappa}
\def\G{\Gamma}
\def\E{\Bbb E}
\def\P{\Bbb P}
\def\k{\kappa}
\def\tG{\tilde{\G}}
\document
\topmatter
\title expected number of distinct  part sizes   in a random integer 
composition 
\endtitle
\leftheadtext{}
\rightheadtext{}
\author pawe{\l} hitczenko and   gilbert  stengle\endauthor

\thanks Part of the research of the first  author was carried out while he was
visiting Department of Mathematics of Lehigh University. He would like to thank the Department for its hospitality.\endthanks
\affil  North Carolina State University and Lehigh University \endaffil \address
Department of Mathematics, NC State University, Raleigh, NC 27695-8205,
USA\endaddress \email pawel\@math.ncsu.edu\endemail

\address Department of Mathematics, Lehigh University,  Bethlehem, PA
18015, USA\endaddress \email gas0\@lehigh.edu\endemail
\endtopmatter

\heading 1. Introduction\endheading
In this  note we obtain  precise  asymptotics, as $n\to\infty$, for the
expected number of distinct 
 part sizes in a random composition  of an integer $n$. Let us recall  that
a multiset $\lambda=
\{\lambda_1,\dots,\lambda_k\}$ is a partition of an integer $n$ 
if the $\lambda_j$ are positive integers,  called parts,  such that
$\sum\lambda_j=n$. The values of $\lambda_j$'s are called part sizes.
Compositions are   partitions in which the order of parts is significant. 
Thus,
for example, the  integer 3 admits three partitions, $\{1,1,1\}$, $\{2,1\}$ 
and
$\{3\}$, and four compositions, namely $(1,1,1)$, $(1,2)$, $(2,1)$ and $(3)$
According to our terminology $(1,2)$ is a composition of $3$ in two parts with
sizes $1$ and $2$. In analogy with random partitions, by a random composition
of an integer $n$ we mean a composition of $n$ that is chosen uniformly at
random out of the set of all $2^{n-1}$ compositions of an integer $n$. More
formally, one considers the probability space consisting of the set $C(n)$ of
all compositions of $n$ equipped  with the uniform probability measure. In
this setting, the number of distinct part sizes (or other characteristics) becomes
a random variable whose probabilistic behavior is to be studied.  
 
Investigation of random partitions from this probabilistic perspective
originated with a paper by Erd\"os and Lehner \cite{5} who studied the limiting
distribution of the total number of parts in a random  partition.
Subsequently,  Wilf \cite{11}  found an asymptotic formula for the expected
number of distinct part sizes. Goh
and Schmutz  \cite{7} obtained more precise information on the distribution of
the number of distinct part sizes, namely they established   the central limit
theorem. Recently Corteel, Pittel,
Savage and Wilf \cite{3} obtained a refined version of Wilf's result 
concerning
the expectation of the number of distinct part sizes in a random partition.
Their result allows one to obtain as many terms for the asymptotic expansion 
of
this expectation as one wishes. For example, on ``$o(1)$ level'' this
expectation is   $$\frac{\sqrt{6n}}{\pi}+\frac3{\pi^2}-\frac12+o(1).$$

The aim of this note is to obtain an asymptotics for the same quantity in
the case of random compositions. In order to state our 
result we need some notation. For 
an integer $n$   consider the set $C(n)$ of all
compositions $\k$ of $n$ equipped with the uniform probability measure $\P_n=\P$
(that is, $\P(\k)=2^{-n+1}$ for every $\k\in C(n)$.) For a
composition $\k=(\g_1,\dots,\g_k)$ the number of distinct part sizes,
$D_n(\k)$ is defined by the formula
$$D_n(\k)=1+\sum_{i=2}^kI_{\{\g_i\ne\g_j,\ j=1,\dots,i-1\}},$$
where $I_A$ is the indicator function of the set $A$. We denote the
integration with respect to $\P$ on $C(n)$ by $\E$. We have:
\proclaim{Theorem} As $n\to\infty$, 
$$\E D_n=\log_2n+\frac{\g}{\ln2}-\frac32+g(\log_2 n)+o(1),$$
where $\g$ is Euler's constant and $g$ is mean--zero function of
period 1 satisfying $|g|\le 0.0000016$.
\endproclaim
Thus, the expected
number  of distinct part sizes in a composition of an integer $n$ 
asymptotically behaves
like $\log_2n$ plus a constant plus a small but periodic oscillation. This
oscillatory behavior, which just a few years ago was considered surprising (to
say the least) is by now a well documented and acknowledged feature of
sequences of geometric random variables, see e.g. \cite{2}, \cite{4}, 
\cite{8} \cite{10}.

We wish to observe that the asymptotics for the expected number of
distinct part sizes is the same as the expected length of the longest
run of heads in $n$ tosses of a fair coin, see e.g. \cite{2} or
\cite{8}. Since the size of the largest part is one plus the longest
run of heads it follows that on average one expects to see parts of
all but one sizes between 1 and the largest size (or runs of heads of
all but one lengths between 1 and the longest run). 

\heading 2. outline of a proof\endheading
Quite often results like this are obtained
through careful analysis of the the generating function. We will use a
different approach. We will view random composition as (essentially) randomly
stopped sequence of i.i.d. geometric random variables and we will express the
number of distinct part sizes as a function of this sequence.   This will
allow for  direct and straightforward estimates. The same
approach was used successfully in \cite{9} to handle a problem in which
generating function approach was apparently futile. We believe that this
technique will prove useful in many other problems concerning random
compositions.  Our proof in a natural way splits in the two steps. In
the first we will use the afore--mentioned representation and
probabilistic estimates to extract the
main contribution to $\E D_n$. Namely, we have
\proclaim{Proposition 1} As $n\to\infty$, 
$$\E
D_n=\sum_{m=1}^\infty\Big\{1-\big(1-\frac1{2^m}\big)^{\frac{n}2(1+o(1)}\Big\}+o(1).$$
\endproclaim

The second, purely analytical step is to analyse the asymptotic
behavior 
of the infinite sum above. This goal could be accomplished by applying
the so -- called Rice method (see e.g. \cite{6} for a very good
description and examples). Since this method requires some tools from
complex analysis we decided to take a different route. As a result
our analysis is completely elementary (thus, making this paper fully
accessible to  advanced undergraduates, for example.) To facilitate our analysis we define 
$$f(x)=
\sum_{m=1}^\infty\Big\{1-\big(1-\frac1{2^m}\big)^{2^x}\Big\}.$$
With this definition we will show that  $f(x)$ tends to a limit as
$x$ tends to infinity along sequences of the form $\{x_0 + k\}_{k\in
Z}$, but does not possess a unrestricted limit as
$x\rightarrow\infty.$ More specifically, we have 
\proclaim{Proposition 2} For large positive $k$
$$f(x+k) = x + k +\gamma/\ln2-1/2 + g(x) + o(2^{-x-k})$$
where $\gamma$ is Euler's constant and
$$ g(x)= -x - \gamma/\ln2+ 1/2 -\sum^0_{m=-\infty }\exp(-2^{-m+x}) +\sum^\infty_{m=1}(1-\exp(-2^{-m+x}))$$
is a nonconstant, zero-mean function of period 1 satisfying $\mid g(x)\mid \le .0000016.$
\endproclaim
Clearly, Theorem follows by combining Propositions 1 and 2.

\heading 3. proof of proposition 1\endheading

Central to our approach is the following proposition 
\proclaim{Proposition 3} Let
 $\G_1,\G_2\dots$
be i.i.d. geometric random variables with parameter $1/2$  (that is
$\P(\G_1=j)=2^{-j}$, $j=1,2\dots$) and define  $$\tau=\inf\{k\ge1:\
\G_1+\G_2+\dots+\G_k\ge  n\}.$$ Then, the distribution of a randomly 
chosen composition in $C(n)$ is  
given                                                                      
by $$(\G_1,\G_2,\dots,\G_{\tau-1},n-\sum_{j=1}^{\tau-1}\G_j).$$
\endproclaim 
This proposition is nothing more than a reiteration of a known (see
e.g. 
\cite{1}) 
connection between compositions of integers and $\{0,1\}$ -- valued sequences. Namely, a composition $\k=(\g_1,\dots,\g_k)$ of an integer $n$ into
parts $\g_1,\dots,\g_k$ is associated with a string of 0's and 1's of length
$n$  as follows: there is a 1 on the $n$th place and the  numbers
$\g_1,\dots,\g_k$  are ``waiting times'' for the first, second,$\dots,$ and
$k$th appearance of  1. (For example, the composition $(1,2,3,1,1)$ of 8
corresponds to the string $10100111$ while $(4,2,2)$ corresponds to
$00010101$.)  Choosing a composition at random amounts to having  the  $0$'s and $1$'s on the first $n-1$ 
places
occur according to a binomial   $\text{Bin}(n-1,1/2)$ law. We refer to \cite{9}
for more details.  
 Let $\tG_i(\k)$ denote parts of a
randomly chosen composition $\k$, i.e.
$$\tG_i(\k)=\G_i(\k),\quad\text{for}\quad i<\tau(\k)\quad \text{and}\quad  
\tG_{\tau(\k)}(\k)=n-\sum_{i=1}^{\tau(\k)-1}\G_i(\k).$$ 
Note that $\tG_\tau\le\G_\tau$. The 
 expected value of $D_n$  is computed as follows:
 $$\align
\E D_n  &=1+\E \sum_{j=2}^\tau I_{\tG_i\ne\tG_j;\ j=1,\dots,i-1}\\&=1+
\E \sum_{j=2}^{\tau-1} I_{\G_i\ne\G_j;\ j=1,\dots,i-1}+
\P(\tG_{\tau}\ne\G_1,\dots,\G_{\tau-1}).\endalign$$
We will first show that the last probability is negligible. This is because
$\tau$ being a $1+\text{Bin}(n-1,1/2)$ random variable satisfies the bound 
$$\P(|\tau-\E\tau|\ge t)\le 2\exp(-2t^2/(n-1)),$$
so that 
$$\P(|\tau-\E\tau|\ge \sqrt{(n-1)\log n})\le 2\exp\{-2\frac{(n-1)\log
n}{n-1}\}=O(1/n^2).$$
Therefore, letting $t_n=\sqrt{(n-1)\log n}$ and then  
$$n_0\sim\E\tau-t_n=(n+1)/2-\sqrt{(n-1)\log n},$$
and $$n_1\sim
\E\tau+t_n=(n+1)/2+\sqrt{(n-1)\log n},$$ we get 
$$\align\P(\tG_{\tau}\ne\G_1,&\dots,\G_{\tau-1})
\le\sum_{j=1}^n
\P(\tG_\tau=j,\ \G_1,\dots,\G_{\tau-1}\ne j)
\\&\le\sum_{j=1}^n
\P(\tG_\tau\ge j,\ \G_1,\dots,\G_{\tau-1}\ne j)\\&
\le\sum_{j=1}^n
\P(\G_\tau\ge j,\ \G_1,\dots,\G_{\tau-1}\ne j)
\\&\le\sum_{j=1}^n \P(\G_\tau\ge j,\ \G_1,\dots,\G_{\tau-1}\ne j,\
|\tau-\E\tau|\le t_n) +n\P(|\tau-\E\tau|\ge t_n)\\&
\le\sum_{j=1}^n\sum_{k=n_0}^{n_1}
\P(\G_k\ge j,\G_1,\dots,\G_{k-1}\ne j,\ \tau=k)+O(1/n)
\\&\le\sum_{k=n_0}^{n_1}\sum_{j=1}^\infty
\P(\G_k\ge j, \G_1,\dots,\G_{k-1}\ne j)+O(1/n)
\\&=\sum_{k=n_0}^{n_1}\sum_{j=1}^\infty\frac1{2^{j-1}}
\big(1-\frac1{2^j}\big)^{k-1}+O(1/n)
\\&\le\sum_{k=n_0}^{n_1}
C\int_0^\infty\frac1{2^x}\big(1-\frac1{2^{x+1}}\big)^{k-1}dx+O(1/n)
\\&\le\sum_{k=n_0}^{n_1} \frac{C}{k}+O(1/n)\le C\frac{\sqrt{n\log n}}n\to0
,\endalign$$
as $n\to\infty$.  
 As for the other term, we have  $$\align
\E\sum_{i=2}^{\tau-1}I_{\G_i\ne \G_j,\ j<i}&\le 
\E\big(\sum_{i=2}^{\tau-1}I_{\G_i\ne \G_j,\ j<i}\big)I_{|\tau-\E\tau|\le
t_n}\\&\quad +\E\big(\sum_{i=2}^{\tau-1}I_{\G_i\ne \G_j,\
j<i}\big)I_{|\tau-\E\tau|> t_n} \\&\le \E\big(\sum_{i=2}^{\tau-1}I_{\G_i\ne
\G_j,\ j<i}\big)I_{|\tau-\E\tau|\le t_n} +(n-2)\P(|\tau-\E\tau|\ge 
t_n).\endalign$$
The second term is bounded above by $C/n$  and, of course,  tends to 0 as 
$n\to
\infty$. For the first one we have:  
$$\align
\E\big(\sum_{i=2}^{\tau-1}I_{\G_i\ne
\G_j,\ j<i}\big)&I_{|\tau-\E\tau|\le t_n}\le
\E\big(\sum_{i=2}^{n_1}I_{\G_i\ne \G_j,\ j<i}\big)I_{|\tau-\E\tau|\le t_n}
\\&\le\E\sum_{i=2}^{n_1}I_{\G_i\ne
\G_j,\ j<i}=\sum_{i=2}^{n_1}\P(\G_i\ne
\G_j,\ j<i)
\endalign$$
Similarly, 
$$
\align
\E\big(\sum_{i=2}^{\tau-1}&I_{\G_i\ne
\G_j,\ j<i}\big)I_{|\tau-\E\tau|\le t_n}\ge
\E\big(\sum_{i=2}^{n_0}I_{\G_i\ne \G_j,\ j<i}\big)I_{|\tau-\E\tau|\le t_n}\\&
=\E\sum_{i=2}^{n_0}I_{\G_i\ne
\G_j,\ j<i}-\E\sum_{i=2}^{n_0}\E\big(\G_i\ne
\G_j,\ j<i\big)I_{|\tau-\E\tau|> t_n}\\&\ge \sum_{i=2}^{n_0}\P(\G_i\ne
\G_j,\ j<i)-n\P(|\tau-\E\tau|\ge t_n)\\&= \sum_{i=2}^{n_0}\P(\G_i\ne
\G_j,\ j<i)-O(n^{-1})
\endalign$$

We will now fix  $k$ and approximate $\sum_{i=2}^k\P(\G_i\ne
\G_j,\ j<i)$ as follows
$$
\align
\P(\G_i\ne\G_j,\ j<i)=\sum_{m=1}^\infty\P(\G_i=m,\ \G_j\ne m;\ j<i)=
\sum_{m=1}^\infty\frac1{2^m}\Big(1-\frac1{2^m}\Big)^{i-1}.\endalign
$$
Hence, by summing up over $i$ we get:
$$
\align
\sum_{i=2}^k\sum_{m=1}^\infty&\frac1{2^m}\Big(1-\frac1{2^m}\Big)^{i-1}=
\sum_{m=1}^\infty\frac1{2^m}\sum_{i=1}^{k-1}\Big(1-\frac1{2^m}\Big)^i\\&=
\sum_{m=1}^\infty\frac1{2^m}\Big(1-\frac1{2^m}\Big)\frac{1-(1-2^{-m})^{k-1}}
{1-(1-2^{-m})}\\&=\sum_{m=1}^\infty\Big(1-\frac1{2^m}\Big)
\left(1-\Big(1-\frac1{2^m}\Big)^{k-1}\right)
=\sum_{m=1}^\infty\Big\{\big(1-\frac1{2^m}\big)-\big(1-\frac1{2^m}\big)^k
\Big\}\\&=\sum_{m=1}^\infty\Big\{1-\big(1-\frac1{2^m}\big)^k\Big\}-\sum_{m=1}^\infty\frac1{2^m}=\sum_{m=1}^\infty\Big\{1-\big(1-\frac1{2^m}\big)^k\Big\}-1.\endalign
$$
It follows that after ignoring terms of order $o(1)$ we have 
$$\sum_{m=1}^\infty\Big\{1-\big(1-\frac1{2^m}\big)^{n_0}\Big\}\le \E D_n
\le \sum_{m=1}^\infty\Big\{1-\big(1-\frac1{2^m}\big)^{n_1}\Big\}.$$
since both $n_0$ and $n_1$ are of the form $\frac{n}2(1+o(1))$
Proposition 1 follows.

\heading 4. proof of proposition 2\endheading 
 Recall that 
 $$f(x)=\sum_{m=1}^\infty\Big\{1-\big(1-\frac1{2^m}\big)^{2^x}\Big\}.$$ 
Throughout the proof we can suppose that $0\le x <1.$ We first give a
simple argument which gives the limiting behavior of $f(x)$ without, however, yielding an estimate for the rate of convergence. We re-index the sum by $m+k$ to obtain
$$f(k+x)-k-x = -x -\sum^0_{m=-k+1}(1-2^{-m-k})^{2^{k+x}}+\sum^\infty_{m=1}(1-(1
-2^{-m-k})^{2^{k+x}})       .$$
Permuting summation and limits as $k \rightarrow \infty $ yields 
$$f(x+k)- k -x
 = -x-\sum^0_{m=-\infty }\exp(-2^{-m+x}) +\sum^\infty_{m=1}(1-\exp(-2^{-m+x}))+o
(1).$$
But this step is justified by dominated convergence using the
following majorizing convergent series of positive terms independent
of $k$:  
 $$\sum^0_{m=-k+1}(1-2
^{-m-k})^{2^{k+x}} \ll \sum^0_{m=-\infty}\exp({-2^{-m+x}}) $$
and
$$\sum^\infty_{m=1}(1-(1-2^{-m-k})^{2^{k+x}}) \ll \sum^\infty_{m=1}(1-\exp({-2^{-m+x+1}}))  .$$
These follow from the  estimates    $\exp(-2ab) \le (1-b/\lambda)^{a\lambda}\le
\exp(-ab)$ if $\lambda a>0$ and $b/\lambda \le 1/2$ with $\lambda = 2^k, a=2^x$
and $b=2^{-m}.$

The series thus established as the limit of $f(k+x)-(k+x)$  defines a function of period 1. Denoting its mean by $c$ and its zero--mean part by $g(x)$ we have
$$f(x+k) = x+ k + c +g(x) +o(1) = x+k+c + g(x+k)+o(1).$$

To obtain the finer estimate stated in the proposition we use the higher order estimate
$\exp(-ab-ab^2/\lambda) \le (1-b/\lambda)^{a\lambda}$ if $\lambda a>0$
and 
$b/\lambda \le 1/2.$
We must bound
$$\align f(x+k)-(x+c+g(x))& =  \sum^{-k}_{m=-\infty} \exp({-2^{-m+x}})
\\&\qquad + \sum^{\infty}_{m=-k+1}\big\{ \exp({-2^{-m+x}}) 
-(1-2^{-m-k})^{2^{k+x}}\big\}.
\endalign$$
The first sum can be rewritten as
$$\sum^{\infty}_{m=0} \exp({-2^{m+k+x}})=\exp{(-2^{k+x})}\sum^{\infty}_{m=0} \exp{  ( -(2^m-1)2^{k+x}  )   }$$
which is bounded by
$$\exp{(-2^{k+x})}\sum^{\infty}_{m=0} \exp{(-(2^m-1))}$$
and thus makes an exponentially small contribution to an error term of $O(2^{-k-x}).$

The second sum consists of positive terms and is bounded above by
$$\sum^{\infty}_{m=-\infty}
\{\exp{(-2^{-m+x})}-\exp{(-2^{-m+x}-2^{-2m+x-k})}\}.$$
Then the inequality $\exp{(-a)}-\exp{(-a-b)} \le b\exp{(-a)}$ for positive $a$ and $b$ gives the bound
$$2^{-x-k}\sum^{\infty}_{m=-\infty}2^{-2m+2x}
\exp{(-2^{-j+x})}.$$ This bound has the form $2^{-x-k}h(x)$ where $h$
is a 
periodic function of $x$ and is therefore bounded by a constant. This establishes the asserted rate of convergence.

It remains to calculate the mean of $g$. The mean of $-x$ is $-1/2$
and the mean of the residual series is
$$c_0= -\sum^0_{m=-\infty }\int^1_0 \exp(-2^{-m+x})dx
+\sum^\infty_{m=1}
\int^1_0(1-\exp(-2^{-m+x}))dx.
$$
On the $m$--th summand the change of variable $u=2^{-m+x}$ gives
$$c_0\ln2=  -\sum^0_{m=-\infty }\ \int^{2^{-m+1}}_{2^{-m}}
\exp(-u)\frac{du}{u} +\sum^\infty_{m=1}\int^{2^{-m+1}}_{2^{-m}}(1-\exp(-u)) \frac
{du}{u}
$$
or
$$c_0\ln2=  -\int^{\infty}_{1}
\exp(-u)\frac{du}{u} +\int^{1}_{0}(1-\exp(-u)) \frac{du}{u}.
$$
Integrating each integral by parts yields a single integral
$$-\int^{\infty}_{0} \exp(-u)\ln udu $$
which is a well--known integral representing Euler's constant.

We remark that a little more similar reasoning shows that the periodic function
$$h(x) = \sum^{\infty}_{m=-\infty}2^{-2m+2x}
\exp{(-2^{-m+x})}$$
appearing in the rate estimate overestimates the error asymptotically by a factor of two and, in fact,
$$f(x+k) = x + k +\gamma/\ln2-1/2 + g(x) + h(x)2^{-x-k-1}+ O(2^{-2x-2k}).$$

The bound on $g$ and its nonconstant character are easily checked numerically although calculations dealing with $f$ rather than analytically derived asymptotic
 forms are rather sensitive. Alternately its complex Fourier coefficients are easily obtained by a simple variant of the calculation of $c_0$ and have the form
$c_k = \Gamma(2\pi ki/\ln2)/\ln2.$  These are nonzero, small and decrease geometrically in magnitude. For example
$2\mid c_1\mid = .00000157316$ accounts
 for the maximum contribution of the first harmonic while for the second harmonic
$2\mid c_2\mid < 10^{-12}.$

\enddocument